\RequirePackage{thm-restate, amsmath}
\documentclass[a4paper]{article}
\usepackage[utf8]{inputenc}
\usepackage[T1]{fontenc}
\usepackage[english]{babel}
\usepackage{authblk}
\usepackage{amssymb}
\usepackage{xspace}
\usepackage{hyperref}
\usepackage{verbatim}
\usepackage{tikz}
\usepackage[procnumbered,linesnumbered,ruled,vlined]{algorithm2e}
\usepackage{tabularray}
\usepackage{cleveref}
\usepackage[textsize=footnotesize,color=green!40]{todonotes}
\usepackage{stmaryrd}
\usepackage{comment}
\usepackage{apxproof}
\usepackage[left= 3 cm,right=3 cm,top=3 cm,bottom=3 cm]{geometry}
\usepackage{subcaption}

\tikzstyle{noeud}=[circle,inner sep=2, minimum size =3 pt, line width = 1pt, draw=black, fill=white]

\newcommand{\outcomeP}{$\mathcal{P}$\xspace}
\newcommand{\outcomeN}{$\mathcal{N}$\xspace}
\newcommand{\gr}{\mathcal{G}}
\DeclareMathOperator{\CSG}{CSG}
\DeclareMathOperator{\SUB}{SUB}
\DeclareMathOperator{\mex}{mex}
\DeclareMathOperator{\opt}{opt}
\newcommand{\csgset}[1]{$\CSG(\{#1\})$\xspace}
\newcommand{\csgS}[1]{$\CSG(#1)$\xspace}
\newcommand{\subgameset}[1]{$\SUB(\{#1\})$\xspace}
\newcommand{\arck}{\textsc{Arc-Kayles}\xspace}
\newcommand{\ndarck}{\textsc{Non-Disconnecting Arc-Kayles}\xspace}
\newcommand{\nodek}{\textsc{Node-Kayles}\xspace}
\newcommand{\ndnodek}{\textsc{Non-Disconnecting Node-Kayles}\xspace}
\newcommand{\cram}{\textsc{Cram}\xspace}
\newcommand{\strat}{\mathcal{S}}

\newtheorem{theorem}{Theorem}
\newtheorem{definition}[theorem]{Definition}

\newtheorem{lemma}[theorem]{Lemma}
\newtheorem{reduction}{Reduction rule}

\title{Complexity and algorithms for \arck and \ndarck\thanks{The research of the second author was supported by the International Research Center "Innovation Transportation and Production Systems" of the I-SITE CAP 20-25 and by the ANR project GRALMECO (ANR-21-CE48-0004). The research of the third author was supported by the ANR project P-GASE (ANR-21-CE48-0001-01) and the Kempe Foundation Grant No. JCSMK24-515 (Sweden).} \thanks{A short version of this paper will appear at WALCOM 2026~\cite{burke2026complexity}.}}
\author[1]{Kyle Burke}
\author[2,3]{Antoine Dailly}
\author[4]{Nacim Oijid}
\affil[1]{Florida Southern College, Lakeland, FL 33801, USA}
\affil[2]{Université Clermont-Auvergne, CNRS, Mines de Saint-Étienne, Clermont-Auvergne-INP, LIMOS, 63000 Clermont-Ferrand, France}
\affil[3]{Université Clermont Auvergne, INRAE, UR TSCF, 63000, Clermont-Ferrand, France}
\affil[4]{Department of Mathematics and Mathematical Statistics, Umeå University, Sweden}

\date{}

\begin{document}

\maketitle

\begin{abstract}
	\arck is a game where two players alternate removing two adjacent vertices until no move is left, the winner being the player who played the last move. Introduced in 1978, its computational complexity is still open. More recently, subtraction games, where the players cannot disconnect the graph while removing vertices, were introduced. In particular, \textsc{Arc-Kayles} admits a non-disconnecting variant that is a subtraction game. We study the computational complexity of subtraction games on graphs, proving that they are PSPACE-complete even on very structured graph classes (split, bipartite of any even girth). 
	We give a quadratic kernel for \textsc{Non-Disconnecting Arc-Kayles} when parameterized by the feedback edge number, as well as polynomial-time algorithms for clique trees and a subclass of threshold graphs.  We also show that a sufficient condition for a second player-win on \textsc{Arc-Kayles} is equivalent to the graph isomorphism problem.
\end{abstract}

%Introduction
\section{Introduction}
\label{sec-introduction}

\arck is one of the many games introduced by Schaefer in his seminal paper on the computational complexity of games~\cite{schaefer1978complexity}. It is a two-player, information-perfect, finite vertex deletion game, in which the players alternate removing two adjacent vertices and all their incident edges from a graph. The game ends when the graph is either empty or reduced to an independent graph, \textit{i.e.}, a graph containing no edges, and the first player unable to play loses. Since the game is fully deterministic, we consider that both players play perfectly.

While Schaefer proved the PSPACE-completeness of many games, the complexity of \arck surprisingly remains open. Some FPT and XP algorithms have been proposed, and the game has been studied on some specific graph classes~\cite{huggan2016polynomial}, but the results are quite scarce, the game is still generally open even on subdivided stars with three paths.

\arck can also be seen as an \emph{octal game} played on a graph. Octal games are taking-breaking games played on heaps of counters, where the octal code defines how many counters one can take from a heap, as well as how the heap can be left after the move (empty, non-empty, split in two non-empty heaps). An overview of octal games can be found in volume~1 of~\cite{berlekamp2004winning}. In~\cite{beaudou2018octal}, octal games were expanded to be played on graphs, where removing counters means removing a connected subgraph, and splitting a heap means disconnecting the graph. Under this definition, \arck is the octal game $\mathbf{0.07}$.

Non-disconnecting octal games, also called \emph{subtraction games}, have also been studied on graphs~\cite{dailly2019connected}: for a non-empty set of integers $S$, the \emph{connected subtraction game on $S$}, \csgS{S}, is played by removing connected subgraphs of order $k$ (where $k \in S$) without disconnecting the graph. The non-disconnecting condition of subtraction games allows for better results, compared to \arck, on constrained classes such as trees. Note that \ndarck (for \textsc{Non-Disconnecting Arc-Kayles}) is \csgset{2}.

The aim of this paper is to further the study of the computational complexity of subtraction games, particularly \ndarck. We present several PSPACE-completeness results for subtraction games, a sufficient condition for symmetry strategies for \arck, and tractability results for \ndarck: an FPT algorithm parameterized by the feedback edge number, \textit{i.e.}, the minimum size of a set $S$ of edges such that $G - S$ contains no cycle. as well as several polynomial-time algorithms on structured graph classes.

\medskip
\noindent\textbf{Combinatorial Game Theory terminology.}
\arck and subtraction games are impartial games, a subset of combinatorial games which have been extensively studied. Since both players have the same possible moves, impartial games can have two possible \emph{outcomes}: either the first player has a winning strategy (in which case the game has outcome \outcomeN) or the second player does (in which case the game has outcome \outcomeP). A position $G'$ that can be reached by a move from a position $G$ is called an \emph{option} of $G$, the set of options of $G$ is denoted by $\opt(G)$. If every option of $G$ has outcome \outcomeN, then $G$ has outcome \outcomeP; conversely, if any option has outcome \outcomeP, then $G$ has outcome \outcomeN. A refinement of outcomes are the \emph{Sprague-Grundy values}~\cite{grundy1939mathematics,sprague1935mathematische}, denoted by $\gr(G)$ for a given position $G$. The value can be computed inductively, with $\gr(G)=\mex(\opt(G))$, where $\mex(S)$ is the smallest non-negative integer not in $S$. A position $G$ has outcome \outcomeP if and only if $\gr(G) = 0$. Furthermore, \arck and connected subtraction games are \emph{vertex deletion games}, a family of games where the players alternate removing vertices from a graph following some fixed constraints.
For more information and history about combinatorial games, we refer the reader to~\cite{albert2019lessons,berlekamp2004winning,conway2000numbers,siegel2013combinatorial}.

\medskip
\noindent\textbf{\arck: a difficult game.}
Introduced in 1978~\cite{schaefer1978complexity}, \arck's complexity remains open to this day. However, there have been some algorithmic results.
Determining the winner of \arck is FPT when parameterized by the number of rounds~\cite{lampis2014computational} and the vertex-cover number, and XP when parameterized by neighborhood diversity~\cite{hanaka2023winner}. \arck on paths is exactly the octal game $\mathbf{0.07}$, and its Grundy sequence is well-known to have period~34 and a preperiod of~68~\cite{guy1956g}. \arck on grid graphs is exactly the game \cram~\cite{gardner1974mathematical}, and thus there are winning strategies for the first player on even-by-odd grids and for the second player on even-by-even grids (by an easy symmetry argument). The game remains open on odd-by-odd grids; Sprague-Grundy values for grids with at most~30 squares~\cite{uiterwijk2018construction} and of size $3 \times n$ for $n$ up to 21 and a few other odd-by-odd grids~\cite{lemoinecomputation,uiterwijk2019solving} have been computed without any regularity appearing. Finally, the game seems to be hard to solve even on trees: an $\mathcal{O}(2^{\frac{n}{2}})$ enumerative algorithm was proposed~\cite{hanaka2023winner}, later improved to $\mathcal{O}(1.3831^n)$~\cite{hanaka2024faster}. While the game remains open on subdivided stars with three paths, one of which has size~1~\cite{huggan2016polynomial}---by fixing the length of the second path and having the third one change---the authors conjectured that the sequence of Sprague-Grundy values thus obtained would ultimately be periodic. However, the Sprague-Grundy values of \arck are unbounded in the general case~\cite{dailly2019generalization}.

\medskip
\noindent\textbf{Subtraction games on graphs.}
Subtraction games have been extensively studied on heaps of counters, and there are several regularity results on them~\cite{albert2019lessons,berlekamp2004winning}.  It is important to note that these are different rulesets than \textsc{NimG} or \textsc{GraphNim}, games where \textsc{Nim}-heaps are embedded on the nodes or vertices of the graph~\cite{DBLP:journals/tcs/Fukuyama03,DBLP:journals/tcs/Fukuyama03a,DBLP:journals/corr/abs-1101-1507}. Subtraction games on graphs were introduced in~\cite{dailly2019connected}.
They are of particular interest since they restrict the types of moves available on sparser classes (such as trees), and because they cannot be solved through a simple symmetry strategy started by cutting the graph in two isomorphic subgraphs.  Hence, their analysis differs fundamentally from the analysis of vertex deletion games that allow disconnecting.

In~\cite{dailly2019connected}, the authors proved that, by fixing $G$ and a vertex $u$, the Sprague-Grundy values of the sequence of graphs obtained by appending a path of increasing length to $u$ are ultimately periodic. They also found regularity results for the family \csgset{1,\ldots,N} on subdivided stars with three paths, one of which has size~1, in contrast with \arck, which they used as base cases for subdivided stars with specific values of $N$. In~\cite{beaudou2018octal}, the authors gave polynomial-time algorithms for the game \csgset{1,2} on subdivided stars and bistars.

As for the specific case of \ndarck\footnote{A playable version is available at \url{https://kyleburke.info/DB/combGames/nonDisconnectingArcKayles.html}.}, which is the subtraction game \csgset{2}, it is known to be polynomial-time solvable on trees, wheels and grids of height at most~3~\cite{dailly2018criticalite}.

\medskip

\noindent\textbf{Our results and outline of the paper.}

Our focus here is the complexity of the problem that takes as input a graph and outputs the outcome of the considered game played on this graph.

We further expand algorithmic results on \arck and subtraction games on graphs.
Our first contribution is a general PSPACE-completeness reduction for most subtraction games:

\begin{restatable}{theorem}{csgPSPACEcompleteMAIN}
	\label{thm-csgPSPACEcompleteMAIN}
	If $S$ is finite and $1 \not\in S$, then, \csgS{S} is PSPACE-complete, even on bipartite graphs of any given even girth.
\end{restatable}

We also prove that \csgset{k} is PSPACE-complete on split graphs, a very structured graph class. Note that \ndarck is in this family of games.

\begin{restatable}{theorem}{csgPSPACEcompleteSPLIT}
	\label{thm-csgPSPACEcompleteSPLIT}
	For $k \geq 2$, \csgset{k} is PSPACE-complete, even on split graphs.
\end{restatable}

Both of those results are proved in \Cref{sec-hardness}, where we also show that deciding whether a symmetry strategy can be applied for \arck is as hard as the Graph-Isomorphism problem (GI-hard), even on bipartite graphs.

Finally, in \Cref{sec-polyTimeNDArcKayles}, we give polynomial-time algorithms for deciding the winner of \ndarck on several structured graph classes. Our main result is on parameterized complexity. 
Note that the parameters used for \arck (vertex cover and neighbourhood diversity) are linked to the number of rounds and do not allow for connectivity conditions, and as such would not be a good fit. On the contrary, since \ndarck is trivial on trees, the feedback edge number is a natural parameter:

\begin{restatable}{theorem}{feedbackEdgeNumber}
	\label{thm: ndarck parameterized fen}
	\ndarck has a quadratic kernel parameterized by the feedback edge number of the graph.
\end{restatable}

We also study clique trees and subclasses of threshold graphs.

%Classical complexity
\section{Hardness results for vertex deletion games}
\label{sec-hardness}

Before we get to the hardness results, we point out that \csgset{1} is a strategy-free game.  Every vertex can be removed over the course of a game, so players can count the number of moves remaining.

For our first proof, we need to define \nodek, another impartial vertex deletion game.  Each turn, the current player selects one vertex on the graph that is not adjacent to an already selected vertex (there is no distinction between the players' selected vertices). Just as in \arck, under normal play the last player to move wins. Note that, while \nodek is a vertex deletion game (it can be seen as alternately removing one vertex and all its neighbors), it cannot be expressed as a subtraction or octal game on graphs, since the deletion of vertices is restricted by the graph itself and not by a specific set of integers. Recall that the \emph{girth} of a graph is the length of its smallest cycle.

\csgPSPACEcompleteMAIN*

\begin{proof}
	First, note that \csgS{S} is clearly in PSPACE.
	To prove completeness, we reduce from \nodek, which was shown to be PSPACE-complete in~\cite{schaefer1978complexity}.
	
	Let $G(V,E)$ be a position of \nodek. For a given finite set $S$ such that $1 \not\in S$, we construct a position $G'(V',E')$ of \csgS{S} such that the outcomes of $G$ for \nodek and of $G'$ for \csgS{S} are the same.
	
	Let $M = \max(S)$.
	We construct $G'$ from $G$ as follows:
	\begin{itemize}
		\item Create a so-called \emph{control vertex} $c$, and attach $M+1$ leaves to it;
		\item For every vertex $v \in V$, create a star with center $v'$ and $M-1$ leaves, and add the edge $(cv')$;
		\item For every edge $(uv) \in E$, create a vertex $e_{uv}$, attach $M$ leaves to it, and add the edges $(u'e_{uv})$ and $(v'e_{uv})$.
	\end{itemize}
	The reduction is depicted in \Cref{fig-reduction}. It is clearly polynomial, since we have $|V'|=M|V|+(M+1)|E|+M+2$.
	
	We now prove that playing \nodek on $G$ is exactly the same as playing \csgS{S} on $G'$.
	
	First, note that the only possible moves on $G'$ consist of removing a vertex $v'$ and its attached leaves.
	Indeed, since $1 \not\in S$, it is impossible to play on an induced star of order more than $M$, except if the graph is reduced to exactly a star of order $M+1$. Hence, it is always impossible to remove $c$ or any $e_{uv}$ for $(uv) \in E$, leaving the stars centered on vertices $v'$ as the only possible moves.
	
	Furthermore, if a vertex $u'$ has been removed, then it is impossible to remove any vertex $v'$ such that $(uv) \in E$, since we cannot remove the star centered on $e_{uv}$, as doing so would disconnect the graph.
	
	Hence, the only possible moves for \csgS{S} on $G'$ are on vertices equivalent to those of $G$, and playing on such a vertex (removing the star it is the center of) prevents either player from playing on any of its neighbors. This proves that the moves of \csgS{S} on $G'$ are exactly the ones of \nodek on $G$, and thus the game trees (and thus the outcomes and the Grundy values) of those two games are the same.
	
	Note that $G'$ has girth~4, since for any edge $(uv)$ in $G$ there is a 4-cycle $u'cv'e_{uv}$ in $G'$, and that $G'$ is bipartite, since any odd cycle in $G$ now contains twice as many vertices since each edge of $G$ has been subdivided, and no cycle containing $c$ can be of odd length.
	Furthermore, the girth of $G$ can be increased (but will remain even) by subdividing every edge enough times and adding $M+1$ leaves to the newly created vertices: a cycle using $c$ will have length $4+2i$ if we subdivide the edges incident with $c$ $i$ times, afterwards we can subdivide the edges of any more smaller cycles in $G'$ to increase the girth. This transformation does not affect the proof, since the newly created vertices can never be removed, thus proving that $G'$ can have any even girth.
\end{proof}

\begin{figure}[!h]
	\centering
	\begin{tikzpicture}
		\node (orig) at (0,0) {
			\begin{tikzpicture}
				\node[noeud] (0) at (0,0) {};
				\node[noeud] (1) at (2,0) {};
				\node[noeud] (2) at (1,1.5) {};
				\draw (0)--(1)--(2)--(0);
				\draw (0) node[left] {$u$};
				\draw (1) node[right] {$v$};
				\draw (2) node[above] {$w$};
			\end{tikzpicture}
		};
		
		\node (red) at (5,0) {
			\begin{tikzpicture}
				\node[noeud] (0) at (0,0) {};
				\draw (0) node[left] {$u'$};
				\node[noeud] (1) at (4,0) {};
				\draw (1) node[right] {$v'$};
				\node[noeud] (2) at (2,3) {};
				\draw (2) node[left] {$w'$};
				\node[noeud] (01) at (2,0) {};
				\draw (01) node[above] {$e_{uv}$};
				\node[noeud] (02) at (1,1.5) {};
				\draw (02) node[right] {$e_{uw}$};
				\node[noeud] (12) at (3,1.5) {};
				\draw (12) node[left] {$e_{vw}$};
				\node[noeud] (m) at (2,1) {};
				\draw (m) node[above left] {$c$};
				\draw (0)--(01)--(1)--(12)--(2)--(02)--(0);
				\draw (0)--(m);
				\draw (1)--(m);
				\draw (2)--(m);
				\draw (0)--(-0.125,-0.375);
				\draw (-0.125,-0.375) node[noeud] {};
				\draw (0)--(0.125,-0.375);
				\draw (0.125,-0.375) node[noeud] {};
				\draw (1)--(3.875,-0.375);
				\draw (3.875,-0.375) node[noeud] {};
				\draw (1)--(4.125,-0.375);
				\draw (4.125,-0.375) node[noeud] {};
				\draw (2)--(1.875,3.375);
				\draw (1.875,3.375) node[noeud] {};
				\draw (2)--(2.125,3.375);
				\draw (2.125,3.375) node[noeud] {};
				
				\draw (02)--(0.5,1.4);
				\draw (02)--(0.65,1.65);
				\draw (02)--(0.8,1.9);
				\draw (0.5,1.4) node[noeud] {};
				\draw (0.65,1.65) node[noeud] {};
				\draw (0.8,1.9) node[noeud] {};
				
				\draw (12)--(3.5,1.4);
				\draw (12)--(3.35,1.65);
				\draw (12)--(3.2,1.9);
				\draw (3.5,1.4) node[noeud] {};
				\draw (3.35,1.65) node[noeud] {};
				\draw (3.2,1.9) node[noeud] {};
				
				\draw (01)--(1.75,-0.375);
				\draw (01)--(2,-0.375);
				\draw (01)--(2.25,-0.375);
				\draw (1.75,-0.375) node[noeud] {};
				\draw (2,-0.375) node[noeud] {};
				\draw (2.25,-0.375) node[noeud] {};
				
				\draw (m)--(1.75,0.625);
				\draw (m)--(2.25,0.625);
				\draw (m)--(2.25,1.25);
				\draw (m)--(2.4,1);
				\draw (1.75,0.625) node[noeud] {};
				\draw (2.25,0.625) node[noeud] {};
				\draw (2.25,1.25) node[noeud] {};
				\draw (2.4,1.0) node[noeud] {};
			\end{tikzpicture}
		};
		
		\draw[->, line width = 0.7pt, double] (orig)--(red);
	\end{tikzpicture}
	\caption{An illustration of the reduction from \nodek to \csgS{S} with $S = \{2,3\}$.}
	\label{fig-reduction}
\end{figure}

We remark that the non-disconnecting variant of \nodek\footnote{A playable version is available at \url{http://kyleburke.info/DB/combGames/nonDisconnectingNodeKayles.html}.}, where the selected vertices and their neighbors must not induce a separating set (\emph{i.e.}, a set of vertices that disconnects the graph), is also PSPACE-complete:

\begin{theorem}
	\label{thm-nonDisconnectingNodeKayles}
	\ndnodek is PSPACE-complete.
\end{theorem}

\begin{proof}
	We reduce from \nodek, and consider both \nodek and \ndnodek as vertex deletion games. Let $G(V,E)$ be a position of \nodek. We construct a position $G'$ of \ndnodek by replacing, in $V$, every edge $(uv)$ by the following gadget $X_{uv}$:
	\begin{itemize}
		\item Subdivide $(uv)$ twice, let $e_{uv}^1$ and $e_{uv}^2$ be the two vertices thusly created;
		\item Construct a copy of the complete bipartite graph on six vertices $K_{3,3}$, and connect one of its vertices to both $e_{uv}^1$ and $e_{uv}^2$.
	\end{itemize}
	The edge gadget $X_{uv}$ is depicted in \Cref{fig-ndnodek}.
	Afterwards, connect all the $e_{uv}^1$'s and $e_{uv}^2$'s for every edge $(uv) \in E$.
	
	Playing \nodek in $G$ is strictly equivalent to playing \ndnodek in $G'$. Indeed, note that no player can ever select a vertex in a copy of $K_{3,3}$, as doing so would disconnect the two other vertices in its part of the copy from the rest of $G'$. Furthermore, selecting $u$ in $G'$ prevents both players from selecting $v$ for every $(uv) \in E$ afterwards, as this would disconnect the copy of $K_{3,3}$ in $X_{uv}$ from the rest of $G'$. Finally, it is impossible to select a vertex $e_{uv}^1$ or $e_{uv}^2$ from any gadget $X_{uv}$, as doing so would disconnect the copy of $K_{3,3}$ in the gadget from the rest of $G'$. Hence, every move in $G$ has an equivalent move in $G'$ and conversely, so the two game trees, and thus outcomes and Grundy values, are the same.
\end{proof}

\begin{figure}
	\centering
	\begin{tikzpicture}
		\node[noeud,scale=1.5] (u) at (0,2) {};
		\node[noeud] (uv1) at (1,2) {};
		\node[noeud] (uv2) at (2,2) {};
		\node[noeud,scale=1.5] (v) at (3,2) {};
		\draw (u)to(uv1)to(uv2)to(v);
		
		\foreach \I in {0,1,2} {
			\pgfmathsetmacro{\J}{0.5+\I}
			\node[noeud] (a\I) at (\J,1) {};
			\node[noeud] (b\I) at (\J,0) {};
		}
		\foreach \I in {0,1,2} {
			\foreach \J in {0,1,2} {
				\draw (a\I)to(b\J);
			}
		}
		\draw (uv1)to(a1);
		\draw (uv2)to(a1);
		
		\draw (u) node[left,xshift=-1mm,scale=1.5] {$u$};
		\draw (v) node[right,xshift=1mm,scale=1.5] {$v$};
		\draw (uv1) node[above,yshift=1mm,scale=1.25] {$e_{uv}^1$};
		\draw (uv2) node[above,yshift=1mm,scale=1.25] {$e_{uv}^2$};
	\end{tikzpicture}
	\caption{The edge gadget $X_{uv}$ of the reduction of \Cref{thm-nonDisconnectingNodeKayles}}
	\label{fig-ndnodek}
\end{figure}

For the next reduction, we need to define \textsc{Avoid True}, a game played on positive (without negations) boolean formulas and an assignment of all variables such that the formula is false.  (The variables are initially set to false.)  Each turn the current player selects one false variable and switches the assignment to true.  Players may not switch variables that cause the overall formula to evaluate to true.  Under normal play, the last player to move wins.  \textsc{Avoid True} is known to be PSPACE-complete on formulas in disjunctive normal form~\cite{schaefer1978complexity}.
Furthermore, recall that a split graph is a graph where the vertex set can be partitioned into a clique (so every possible edge within the part exists) and an independent set (so no edge exists within this part).

\csgPSPACEcompleteSPLIT*

\begin{proof}
	The proof is obtained by reduction from {\sc Avoid True} (see \Cref{fig-splitREDUCTION}).
	
	Let $\varphi = \underset{i=1}{\overset{m}{\bigvee}} C_i$ be an {\sc Avoid True} formula in disjunctive normal form over the variables $x_1, \dots, x_n$. We construct a graph $G = (V,E)$ as follows:
	
	\begin{itemize}
		\item For any $1 \le i \le n$, we include $k$ vertices $v_i$ and $v_i^1, \dots, v^{k-1}_i$ in $V$ and all the edges $(v_i v^j_i)$ for $1 \le j \le k-1$ in $E$.
		\item For any $1 \le i, j \le n$ we include the edge $(v_i v_j)$ in $E$.
		\item For any $1 \le j \le m$, we include a vertex $c_j$ in $V$, and if $x_i \in C_j$, we add the edge $(v_i c_j)$ to $E$.
	\end{itemize}
	
	The obtained graph (an example for $k=2$ is depicted in \Cref{fig-splitREDUCTION}) is a split graph, as the $v_i$'s are a clique and all the other vertices are an independent set.%\todo{What is a split graph? Do we need to define it? -Kyle}
	
	First note that the only moves available are of the form $\{v_i, v_i^1, \dots, v^{k-1}_i\}$. Indeed, as the $v_i$'s are a vertex cover, any move has to contain at least one $v_i$, and if this move does not contain any of the $v^j_i$ (for $1 \le j \le k-1)$, then this vertex becomes isolated and the graph is no longer connected. 
	
	We now prove that the first (second resp.) player wins in {\sc Avoid True} if and only if they win in \ndarck.
	
	Suppose that the first (second resp.) player has a winning strategy $\strat$ in {\sc Avoid True}. We define their strategy in \ndarck as follows:
	
	\begin{itemize}
		\item Whenever $\strat$ would take a variable $x_i$, they play the set $\{v_i, v_i^1, \dots, v^{k-1}_i\}$
		\item If their opponent plays a set $\{v_i, v^1_i, \dots, v^{k-1}_i\}$, they consider that the variable $x_i$ has been played against $\strat$ in {\sc Avoid True}.
	\end{itemize}
	
	Following this strategy, as $\strat$ is a winning strategy in {\sc Avoid True}, the second (first resp.) player will be the first to satisfy a clause $C_j$. The moves satisfying $C_j$ will therefore be removing the last vertex $v_i \in C_j$ and will disconnect the graph. Therefore, the first (second resp.) player wins.  
\end{proof}

\begin{figure}
	\centering
	\begin{tikzpicture}
		\node[noeud] (x1) at (0,1) {};
		\node[noeud] (x2) at (1,1.5) {};
		\node[noeud] (x3) at (2,1.5) {};
		\node[noeud] (x4) at (3,1) {};
		
		\node[noeud] (y1) at (0,2) {};
		\node[noeud] (y2) at (1,2.5) {};
		\node[noeud] (y3) at (2,2.5) {};
		\node[noeud] (y4) at (3,2) {};
		
		\node[noeud] (c1) at (0,0) {};
		\node[noeud] (c2) at (2,0) {};
		\node[noeud] (c3) at (3,0) {};
		
		\foreach \I in {1,2,3,4} {
			\draw (x\I)to(y\I);
			\foreach \J in {1,2,3,4} {
				\ifthenelse{\I = \J}{}{\draw (x\I)to(x\J);}
			}
		}
		\foreach \I in {1,2} {\draw (c1)to(x\I);}
		\foreach \I in {4,2} {\draw (c2)to(x\I);}
		\foreach \I in {1,3} {\draw (c3)to(x\I);}
		
		\draw (x1) node[left] {$v_1$};
		\draw (x2) node[above left] {$v_2$};
		\draw (x3) node[above right] {$v_3$};
		\draw (x4) node[right] {$v_4$};
		
		\draw (y1) node[left,yshift=1mm] {$v_1^{1}$};
		\draw (y2) node[left,yshift=1mm] {$v_2^{1}$};
		\draw (y3) node[left,yshift=1mm] {$v_3^{1}$};
		\draw (y4) node[left,yshift=1mm] {$v_4^{1}$};

		\draw (c1) node[left,yshift=-1mm] {$c_1$};
		\draw (c2) node[left,yshift=-1mm] {$c_2$};
		\draw (c3) node[right,yshift=-1mm] {$c_3$};
		
	\end{tikzpicture}
	\caption{An illustration of the reduction from \textsc{Avoid True} to \csgS{k} with $k=2$. The clauses are $C_1=(x_1 \land x_2)$, $C_2=(x_2 \land x_4)$ and $C_3=(x_1 \land x_3)$. The $v_i$'s form a clique, and the other vertices an independent set.}
	\label{fig-splitREDUCTION}
\end{figure}

Returning to \arck, we present a natural and sufficient winning condition for the second player, which is used in the case of grids.  When the graph is highly symmetric, the second player may want to apply a symmetry strategy with respect to a given "point" in the graph.
We show that, in the general case, deciding whether such a strategy is possible is as hard as the famous Graph Isomorphism problem. 
Note that such symmetry strategies are often used in games, and are often GI-hard to compute, see~\cite{bensmail2022largest,duchene2024bipartite}.

\begin{definition}
	An automorphism $f$ of a graph $G$ is said to be {\em edge-disjoint} if for any $e \in E(G)$, we have $e\cap f(e) = \emptyset$.
\end{definition}

\begin{theorem}
	\label{thm-automorphismImpliesSymmetricStrategy}
	Let $G = (V,E)$ be a graph such that there exists an edge-disjoint involutive automorphism $f$ of $G$. Then $G$ has outcome \outcomeP in \arck.
\end{theorem}

\begin{proof}
	Let $G$ be a graph, and let $f$ be an involutive automorphism of $G$ satisfying the hypothesis of the Theorem.
	Consider the following strategy for the second player: each time his opponent plays an edge $e$, he answers by playing $f(e)$. This move is always available as $e \cap f(e) = \emptyset$, and $G \setminus V(e, f(e))$ still satisfies the hypothesis of the theorem.
\end{proof}

Unlike \arck, in \ndarck, there are graphs with such an $f$ where the symmetry strategy does not give a winner, such as the cycle on six vertices: it has this property, but the graph has outcome \outcomeN. Indeed, in \ndarck, the move $f(e)$ can be forbidden as it could disconnect the graph.

\begin{theorem}
	\label{thm-automorphismIsGIhard}
	Verifying that a graph $G$ admits an involutive automorphism $f$ such that, for all $e$, $e \cap f(e) = \emptyset$ is GI-hard, even in bipartite graphs.
\end{theorem}

\begin{proof}[of \Cref{thm-automorphismIsGIhard}]
	We provide a reduction from the graph isomorphism problem. Let $G_1, G_2$ be two graphs. If $|V(G_1)| \neq |V(G_2)|$, let $G' = K_2$, two vertices connected by an edge. Otherwise, denote by $n = |V(G_1)|= |V(G_2)|$, and let $G_1', G_2'$ be the graphs obtained from $G_1$ and $G_2$ obtained by adding a universal vertex $u_i$ and a leaf connected to it $v_i$ in $G_i$, and by subdividing all the edges of it once. Let $G' = G_1' \cup G_2'$. Note that $G'$ is bipartite thanks to the subdivision of its edges. We prove that $G_1$ and $G_2$ are isomorphic if and only if $G'$ admits an edge-disjoint involutive automorphism $f$.
	
	Suppose first that $G'$ admits an edge-disjoint involutive automorphism $f$. As $u_1$ and $u_2$ are the only two vertices of degree $n+1$, they have to be the image of the other by $f$. Then, since $v_1$ and $v_2$ both are leaves connected to $u_1$ and $u_2$ by paths of length $2$, and since $f(u_1) = u_2$, we can suppose without loss of generality that $f_1(v_1) = v_2$. 
	
	Now, consider $\phi$ the restriction to $G_1$ of $f$. If $w \in G_1$, since automorphisms preserve distances, we have that $f(\{w\}$ must be at distance $2$ from $f(u_1)$. Therefore, since $f(u_1) = u_2$, $f(w)$ is at distance $2$ from $u_2$ and therefore corresponds to a vertex in $G_2$. This proves that the restriction $G_1$ of $f$ goes to $G_2$, and thus, since $f$ is injective as an automorphism and since $|V(G_1)|= |V(G_2)|$, $\phi$ is an isomorphism between $G_1$ and $G_2$
	
	Suppose now that $G_1$ and $G_2$ admits an isomorphism $f$. Then consider $f': G'_1 \to G'_2$ defined by $$f'(v) =   \left\{
	\begin{array}{ll}
		u_2 & \mbox{if } v=u_1 \\
		v_2 & \mbox{if } v=v_1 \\
		f(v) & \mbox{if } v\in G_1 \\
		u & \mbox{if $v$ is the middle vertex of an edge $(ab)$ and $u$ is the middle vertex of $\left (f(a)f(b) \right )$} 
	\end{array}
	\right.
	$$
	
	Now consider $\phi : G' \to G'$ defined by
	$$\phi(v) =   \left\{
	\begin{array}{ll}
		f'(v) & \mbox{if } v \in G'_1 \\
		f'^{-1}(v) & \mbox{if } v \in G'_2
	\end{array}
	\right.
	$$
	
	By construction $\phi$ is an involutive automorphism of $G'$, and for any $e \in G'$, we have $e\in G'_1$, if and only if $\phi(e) \in G'_2$. therefore $e\cap \phi(e) = \emptyset$. so this automorphism is edge-disjoint. 
\end{proof}

Put more simply, deciding whether the second player can win a game of \arck by applying a symmetry strategy is GI-hard in the general case.

\section{\ndarck on structured graph families}
\label{sec-polyTimeNDArcKayles}

While there are very few polynomial-time algorithms for \arck on structured graph classes, its non-disconnecting variant seems more approachable. In this section, we expand on the results in~\cite{dailly2018criticalite}, where it was shown that trees are trivial to solve. Indeed, as with \csgset{1}, no strategy is required, as all possible moves on a tree will end up being played, and so the Sprague-Grundy value is either~0 or~1. The result can be computed in linear time by maintaining a representation with levels and labeled leaves ensuring that no edge has to be considered more than once, when it is removed.

First, we will show that this result on trees can be extended to tree-like graphs: graphs with bounded feedback edge number and clique trees. Recall that the \emph{feedback edge set} of a graph is a set of edges that, when removed, turn the graph into a forest; the smallest size of a feedback edge set is the \emph{feedback edge number}. Also recall that a graph is a \emph{clique tree} (or \emph{block graph}) if every biconnected component is a clique.

\feedbackEdgeNumber*

To prove \Cref{thm: ndarck parameterized fen}, we present the following reduction rules. A reduction rule is called {\em safe} if it outputs a graph with the same outcome as the original graph. Moreover, since we may have to consider that a move is played, we consider the decision problem ``Is the outcome of \ndarck played on $G$ $o$?'' for $o \in \{\mathcal{N,P}\}$, to handle changes of the current player in the reduction rules.

\begin{reduction}\label{red1}
	If $v \in G$ has at least three leaves attached to it, remove all of them but three.
\end{reduction}

This reduction rule can easily be applied in $\mathcal{O}(|E|)$ since each edge is considered at most twice.

\begin{lemma}
	Reduction rule~\ref{red1} is safe
\end{lemma}

\begin{proof}
	Let $v$ be a vertex of $G$ and let $\ell_1, \dots, \ell_k$ with $k\ge3$ be the leaves attached to it. The first time a leaf $\ell_i$ is played, it must be through the edge $(\ell_iv)$, which would disconnect the other leaves, {and at least two other leaves will remain}. Thus, playing any of the leaves would not be a legal move. This statement remains true if only three leaves are attached to $v$. Thus, we can remove all but three of them without changing the legal moves, hence the outcome does not change.
\end{proof}

Let $C$ be a smallest connected subgraph of $G$ containing all the two-connected components of $G$, and let $T_1, \dots, T_p$ be the connected components of $G\setminus C$. Note that each $T_i$ is a tree and is connected to $C$ by a single edge. Let $r_1, \dots, r_p$ be the vertices of $T_1, \dots, T_p$ connected to $C$. 

\begin{lemma}\label{lemma max number vertices of degree 3 in C}
	There are at most $2k-2$ vertices of degree at least $3$ in $C.$
\end{lemma}

\begin{proof}
	Note that $C$ has minimum degree $2$, is connected and has a feedback edge set of size $k$. Thus, we have $|E| = n-1 + k$. Then, using the handshaking formula, we have:
	\begin{align*}
		2|E| &= \underset{v \in C}{\sum} deg_C(v) \\
		2(n-1+k) &= 2n + \underset{v \in C}{\sum} (deg_C(v) - 2) \\
		2n+2k-2 &\ge 2n + |\{v\in C | deg_C(v) \ge 3\}| \\
		2k-2 &\ge  |\{v\in C | deg_C(v) \ge 3\}| 
	\end{align*} 
\end{proof}

\begin{reduction}\label{red: move tree}
	Let $e_1 = (ab), e_2 = (cd)$ be two non-adjacent edges of $T_{i_1}$ and $T_{i_2}$ respectively (we may have $i_1 = i_2$), such that the moves on $e_1$ and $e_2$ are legal (we may have that the move on $e_2$ is only legal after the move on $e_1$). We transform $G$ into $G\setminus\{a,b,c,d\}$. 
\end{reduction}

$C$ can be computed in time $\mathcal{O}(|E|^2)$ by simply starting with $C = V(G)$ and removing one by one the vertices of degree $1$ in $C$. Then this reduction rule can be applied in time $\mathcal{O}(|E|^2)$ since each time a pair of leaves is considered, one is removed.

\begin{lemma}
	Reduction rule~\ref{red: move tree} is safe.
\end{lemma}

\begin{proof}
	Let $G$ be the original graph and $G'$ be the graph where one move in $T_{i_1}$ and one move in $T_{i_2}$ have been played for some $1 \le i_1, i_2 \le p$. Denote by $e_1 = (ab)$ and $e_2 = (cd)$ the two edges such that $G' = G\setminus\{a,b,c,d\}$. Note that it is possible that $i_1 = i_2$. We prove that $G$ and $G'$ have the same outcome by induction on $n = |V(G)|$. If $n \le 4$, either the graph is $K_{1,3}$ the star with three leaves, or the graph can be emptied and there is nothing to do. Suppose $n \ge 5$. 
	
	Assume first that $o(G') = \mathcal{N}$, \emph{i.e.} the first player wins in $G'$. Let $f_1$ be a winning first move in $G'$. 
	\begin{itemize}
		\item If $f_1 \in C$, the first player also plays $f_1$ in $G$. Let $f_2$ be the answer of the second player in $G$. 
		\begin{itemize}
			\item If $f_2$ is a possible move in $G'$, by the induction hypothesis, by considering $G''$ the current graph, $G''$ and $G''\setminus\{a,b,c,d\}$ have the same outcome. Thus $o(G'') = \mathcal{N}$.
			\item Otherwise, $f_2 \in \{(ab), (cd)\}$. Thus, the first player can play the other edge in $\{(ab),(cd)\}$ and we are left to the game in $G$' in which the first player has played $f_1$ and thus wins. 
		\end{itemize}
		\item If $f_1 =(a'b') \notin C$, then play $e_1$ instead. By the induction hypothesis, applied on $G\setminus\{a,b\}$ with the edges $e_2$ and $f_1$ (possible as $f_1$ is a legal move in $G'$), and since $G'\setminus \{a',b'\}$ has outcome $\mathcal{P}$, $G\setminus\{a,b\}$ has outcome $\mathcal{P}$, which proves that playing $e_1$ was a winning move for the first player.
	\end{itemize}

	Assume now that $o(G') = \mathcal{P}$. Let $f_1 = (a'b')$ be the first edge played in $G$.
	\begin{itemize}
		\item If $f_1 \in \{e_1, e_2\}$, then the second player answers by playing the second edge in $\{e_1, e_2\}$, and reaches the position $G'$ which is winning for them by hypothesis.
		\item If $f_1 \notin \{e_1, e_2\}$, then by the induction hypothesis applied on $G\setminus\{a',b'\}$, with the edges $e_1, e_2$, $G\setminus\{a',b'\}$ has the same outcome as $G' \setminus \{a'b'\}$, which is winning for the second player by hypothesis.
	\end{itemize}
	This proves that the reduction rule preserves the outcome and thus is safe.
\end{proof}

After applications of Reduction rule~\ref{red: move tree}, at most one move can be played in a tree $T_i$. Up to renaming the trees, suppose it is in $T_1$. Note that we cannot get rid of the last edge on which a move is possible since it may change the strategy in the rest of the graph, and it might be that the outcome doesn't change by removing an edge in a tree (see Figure~\ref{fig: removing one edge does not work}).

\begin{figure}
	\centering
	
	\begin{tikzpicture}[scale=1, every node/.style={circle,inner sep=2, minimum size =3 pt, line width = 1pt, draw=black, fill=white}]
		
		\node (a1) at (-1,2) {};
		\node (a2) at (0,2) {};
		\node (a3) at (1,2) {};
		
		\node (c1) at (0,1) {};
		\node (b1) at (1,1) {};
		\node (d1) at (2,1) {};
		
		\node (l1) at (-2,2) {};
		\node (c2) at (-1,1) {};
		\node (l3) at (-2,0) {};
		\node (b2) at (1,0) {};
		\node (d2) at (2,0) {};
		
		\node (l2) at (-2,1) {};
		\node (u) at (-3,1) {$u$};
		\node (v) at (-4,1) {$v$};
		
		\draw (a1) -- (c1) -- (a2);
		\draw (a3) -- (c1) -- (b1);
		\draw (l1) -- (c2) -- (c1) -- (b2) -- (d2) -- (d1) -- (b1);
		\draw (u) -- (l2) -- (c2) -- (l3);
		\draw[line width=.75mm,color = red] (u) -- (v);
	\end{tikzpicture}
	
	\caption{Example of graph having a pendant edge $(uv)$ such that $o(G) = \mathcal N$ and $o(G \setminus \{u,v\}) = \mathcal N$. }
	\label{fig: removing one edge does not work}
\end{figure}

\begin{reduction}\label{red: replacing trees}
	If $F$ is a forest containing at least $2$ vertices in which no move can be played, and attached to some vertex $u$, we replace $T_u$ the subtree induced by $u$ and this forest by a tree $T'$ of Figure~\ref{fig: replacing trees} such that $o(T_u) = o(T')$ and $o(T_u - \{u\}) = o(T'-\{u\})$ (or both $T_u - \{u\}$ and $T'-\{u\}$ are disconnected). 
	
	If $F$ is a forest containing at least $4$ vertices in which exactly one legal move $e$ can be played, and attached to some vertex $u$, we proceed the same way considering $F-\{e\}$ and we then add a path of length~2 attached to a leaf of~$T'$.
\end{reduction}

An example of application of Reduction rule~\ref{red: replacing trees} is provided in Figure~\ref{fig:reduction rule 3}

\begin{figure}
	
	\centering
	\begin{tikzpicture}[scale=1, every node/.style={circle,inner sep=2, minimum size =3 pt, line width = 1pt, draw=black, fill=white}]
		
		% Gros sommet G (rayon = 2cm)
		\draw[thick] (0,0) circle (1.5cm);
		\node[draw = white] at (0,0) {\huge $G$};
		
		% Sommet u sur l'extrémité droite du cercle
		
		\node (u) at (1.5,0) {$u$};

		% Chemin principal du caterpillar
		\node[] (v1) at (2.5,0) {};
		\node[] (v2) at (3.5,0) {};
		\node[] (v4) at (4.5,0) {};
		
		\draw (u)--(v1)--(v2)--(v4);
		
		% Feuilles attachées vers le haut
		\node[] (l1) at (2.5,1) {};
		\node[] (l2) at (3.5,1) {};
		\node[] (l4) at (4.5,1) {};
		\node[] (l5) at (4.5,-1) {};
		
		\draw (v1)--(l1);
		\draw (v2)--(l2);
		\draw (v4)--(l4);
		\draw (v4)--(l5);
		
		\node[draw = white] (fleche) at (5,0) { \Large $\to$};
		
		\draw[thick] (7,0) circle (1.5cm);
		\node[draw = white] at (7,0) {\huge $G$};
		\node (u) at (8.5,0) {$u$};
		
		\node (a) at (10,.5) {};
		\node (b) at (9,0) {};
		\node (c) at (10,-.5) {};
		
		\draw (a) -- (b) -- (c);
		\draw (u) -- (b);

	\end{tikzpicture}

	\caption{An example of application of Reduction rule~\ref{red: replacing trees}. Here, we have $o(T_u) = \mathcal P$ and $o(T_u-\{u\}) = \mathcal N$.}
	\label{fig:reduction rule 3}
\end{figure}

This reduction rule can be applied in time $\mathcal{O}(|V|)$ since each vertex is in at most one forest and, as mentioned before, the outcome of the game played on a tree can be computed in linear time. If no legal move remains in $F$, we present in Figure~\ref{fig: replacing trees} a set of trees covering all the possible outcomes.

\begin{lemma}\label{weak lemma replacing forest}
	Reduction rule~\ref{red: replacing trees} is safe.
\end{lemma}

In order to prove Lemma~\ref{weak lemma replacing forest}, we prove the following stronger lemma (unlike in Reduction rule~\ref{red: replacing trees}, we take the root into account in the number of vertices):

\begin{lemma}\label{lemma replacing trees}
	Let $G_0$ be a graph, and let $(T,r)$ and $(T',r')$ be two rooted trees of order at least $3$ ($5$ resp.), in which no (exactly one resp.) legal move not containing the root can be made. Let $u \in G$ be a vertex. Suppose $o(T) = o(T')$, $T-\{r\}$ is disconnected if and only if $T' - \{r'\}$ is and $o(T - \{r\}) = o(T' - \{r'\})$ otherwise. Let $G$ and $G'$ be the graphs in which $T$ and $T'$ are attached to $G_0$ by identifying $r$ and $u$ or $r'$ and $u$ respectively. Then $o(G) = o(G')$.
\end{lemma}

\begin{proof}
	
	Suppose that the first (second resp.) player has a winning strategy in $G$. We consider the following strategy in $G'$. While there is a vertex in $G_0$, it is not possible to play any edge containing $u$ without disconnecting the remaining vertices of $G_0$ and the remaining vertices of $T$ or $T'$. Thus, there is a natural bijection between the moves of $G$ and the moves of $G'$, since there is at most one move outside $G_0$ that can be made in the bijection if it exists. When the last move is $G_0$ is made, there are only four possibilities:
	\begin{itemize}
		\item At least three vertices remain in $G_0$. In this case, the game ends as it is not possible to claim $u$ and a vertex of $T$ nor $T'$, since it would disconnect the remaining vertices of $G_0$ with the remaining vertices of $T$ or $T'$, which exists since we supposed that they have order at least $5$.
		\item Exactly two vertices remain in $G_0$. Then, since all the moves have been legal, they are exactly $u$ and a leaf $\ell$ connected to it. Since the last move of $G_0$ has been played, playing $(u\ell)$ is not legal, which means that the resulting graph would be $T - \{r\}$, which would be disconnected. But, by the hypothesis this means that $T' - \{r'\}$ is also disconnected and the game is also finished in $G'$.
		\item Exactly one vertex remains in $G_0$. It must be $u$ since, otherwise, $G$ would be disconnected. The game is then equivalent to the game played on $T$, which has the same outcome as the game played on $T'$
		\item No vertex remains in $G_0$. The game is then equivalent to the game played on $T - \{r\}$, which has the same outcome as the game played on $T' - \{r'\}$. \qedhere
	\end{itemize}
\end{proof}

\begin{figure}
	\centering
	
	\begin{subfigure}[b]{.3 \textwidth}
		
		\centering
		
		\begin{tikzpicture}[scale=1, every node/.style={noeud}]
			\node (a) at (0,.5) {};
			\node (b) at (-1,0) {$u$};
			\node (c) at (0,-.5) {};
			
			\draw (a) -- (b) -- (c);
		\end{tikzpicture} 
		
		\caption{A tree $T$ satisfying $o(T) = \mathcal N$ where removing $u$ disconnects the graph. }
	\end{subfigure} \hfill \begin{subfigure}[b]{.3 \textwidth}
		
		\centering
		
		\begin{tikzpicture}[scale=1, every node/.style={noeud}]
			\node (a) at (0,.5) {};
			\node (b) at (-1,0) {$u$};
			\node (d) at (0,0) {};
			\node (c) at (0,-.5) {};
			
			\draw (a) -- (b) -- (c);
			\draw (d) -- (b);
		\end{tikzpicture}
		
		\caption{A tree $T$ satisfying $o(T) = \mathcal P$ where removing $u$ disconnects the graph.}
	\end{subfigure}  \hfill \begin{subfigure}[b]{.3 \textwidth}
		
		\centering

		\begin{tikzpicture}[scale=1, every node/.style={noeud}]
			\node (a) at (0,.5) {};
			\node (b) at (-1,0) {};
			\node (c) at (0,-.5) {};
			\node (b1) at (-3,0) {$u$};
			\node (b2) at (-2.5,0) {};
			\node (b3) at (-2,0) {};
			\node (b4) at (-1.5,0) {};
			
			\draw (a) -- (b) -- (c);
			\draw (b1) -- (b2) -- (b3) -- (b4) -- (b);
		\end{tikzpicture} 
		
		\caption{A tree $T$ satisfying $o(T) = \mathcal N$ and \\ $o(T-\{u\}) = \mathcal N$. \\ {}}
	\end{subfigure}
	
	\begin{subfigure}[b]{.3 \textwidth}
		
		\centering

		\begin{tikzpicture}[scale=1, every node/.style={noeud}]
			\node (a) at (0,.5) {};
			\node (b) at (-1,0) {};
			\node (c) at (0,-.5) {};
			\node (b4) at (-1.5,0) {$u$};
			
			\draw (a) -- (b) -- (c);
			\draw (b4) -- (b);
		\end{tikzpicture} 
		
		\caption{A tree $T$ satisfying $o(T) = \mathcal P$ and \\ $o(T-\{u\}) = \mathcal N$.}
	\end{subfigure} \hfill \begin{subfigure}[b]{.3 \textwidth}
		
		\centering

		\begin{tikzpicture}[scale=1, every node/.style={noeud}]
			\node (a) at (0,.5) {};
			\node (b) at (-1,0) {};
			\node (c) at (0,-.5) {};
			\node (b3) at (-2,0) {$u$};
			\node (b4) at (-1.5,0) {};
			
			\draw (a) -- (b) -- (c);
			\draw (b3) -- (b4) -- (b);
		\end{tikzpicture}

		\caption{A tree $T$ satisfying $o(T) = \mathcal P$ and \\ $o(T-\{u\}) = \mathcal P$.}
	\end{subfigure} \hfill \begin{subfigure}[b]{.3 \textwidth}
		
		\centering

		\begin{tikzpicture}[scale=1, every node/.style={noeud}]
			\node (a) at (0,.5) {};
			\node (b) at (-1,0) {};
			\node (c) at (0,-.5) {};
			\node (b2) at (-2.5,0) {$u$};
			\node (b3) at (-2,0) {};
			\node (b4) at (-1.5,0) {};
			
			\draw (a) -- (b) -- (c);
			\draw (b2) -- (b3) -- (b4) -- (b);
		\end{tikzpicture} 
		
		\caption{A tree $T$ satisfying $o(T) = \mathcal N$ and \\ $o(T-\{u\}) = \mathcal P$.}
	\end{subfigure}
	\caption{A set of small trees with different outcomes.}
	\label{fig: replacing trees}
\end{figure} 

After an exhaustive use of Reduction rule~\ref{red: replacing trees}, the graph consists in a core $C$ of vertices of degree at least 2 and containing all the 2-connected components of the graph and trees of size at most $10$ attached to it. Recall that, using Lemma~\ref{lemma max number vertices of degree 3 in C}, there are at most $2k-2$ vertices having at least $3$ neighbors in $C$. Thus, we only have to reduce long paths of vertices having $2$ neighbors in $C$. To achieve this, we give types to the vertices of $C$ as follows:

\begin{itemize}
	\item type 0: if nothing is attached to it
	\item type 1: if a leaf is attached to it.
	\item type M: if a tree in which there is a legal move from a leaf attached to it.
	\item type B1: if exactly three leaves are attached to it.
	\item type B2: if any other tree is attached to it. 
\end{itemize}

We also denote by $M$ the set of vertices of type M, and by $B$ the set of vertices of type B1 or B2.
Note that, after exhausting Reduction rule~\ref{red: move tree}, there is at most one vertex of type M.
Let $S = \{v \in C| deg_C(v) \ge 3\} \cup M$. Using Lemma~\ref{lemma max number vertices of degree 3 in C}, we have that $|S| \le 2k-1$.

We now prove that there exists a constant $\alpha$ such that we can reduce any path of vertices of type 0, 1 or $B$ to a path of length at most $\alpha$

Let $P$ be a path of vertices of type 0, 1 or B, and of length at least $3$. Let $a,b$ be the two vertices of $S$ connected to the endpoints of $P$, denoted by $u$ and $v$ respectively. For the rest of the proof, we also consider a new outcome $\mathcal X$ corresponding to ``the move is illegal''. Let $\mathcal O = \{\mathcal N, \mathcal P\}$ the classical set of outcome and $\mathcal O' = \mathcal O \cup \{\mathcal X\}$. We define the {\em extended outcome} $\tilde o(P) = (o_1 \times \varepsilon_1, o_2\times \varepsilon_2, o_3\times \varepsilon_3, o_4\times \varepsilon_4, o_5\times \varepsilon_5, O\times \varepsilon_6)$, with $o_1, o_2, o_3, o_4, o_5 \in \mathcal O'$ the outcome of the subtree  where $a, b, \{a,u\}, \{b,v\}, \{a,b\} $ are removed respectively, $O \subset 2^{\mathcal O \times \mathcal O}$ the set of pairs of outcomes possible after a legal move in $P$ on the subtree rooted in $v$ and the one rooted in $u$, and $\varepsilon_1 \dots, \varepsilon_6 \in \{0,1\}$ the parity of the number in the remaining trees. Note that $\tilde o(P)$ can only take at most $3^5*2^4*2^6 = 248832$ different values.

Let $\mathcal P$ be a finite set of paths (with a forest attached to them) such that  $\tilde o(\mathcal P)$ covers all the possible outcomes of $\tilde o$.

Note that, in practice, these paths are paths in $C$, but are trees in $G$.

\begin{reduction} \label{reduction paths}
	While there is a path of vertices of type $0, 1$ or B of length at least $3$, that is not in $\mathcal{P}$ in $G$, replace it with one of $\mathcal{P}$ with the same extended outcome.
\end{reduction}

This reduction rule can be computed in time $\mathcal{O}(|E|^2)$ as for each of the graphs, the outcome of two trees has to be computed. The following Lemma can be proved using the same idea as Lemma~\ref{weak lemma replacing forest}, but with more case analysis.

\begin{lemma}
	\label{lemma reduction paths safe}
	Reduction rule~\ref{reduction paths} is safe.
\end{lemma}

\begin{proof}
	
	Let $G'$ be the graph obtained from $G$ by replacing a path $P$ by a path $P'$ of $\mathcal P$ having the same extended outcome. Let $a,b$ be the two vertices of $C$ that are connected to $P$. Let $u,v$ be the vertices of $P$ connected to $a$ and $b$ respectively in $G$ and $u',v'$ be the vertices of $P'$ connected to $a$ and $b$ respectively in $G'$
	
	Suppose that the first (second resp.) player has a winning strategy in $G$, we consider the following strategy in $G'$. While the optimal move in $G$ is in $G\setminus (P \cup \{a,b\})$ and while the opponents answers in $G'\setminus (P' \cup \{a,b\})$ , play the same edge in $G'$ or consider he answers the same edge in $G\setminus (P \cup \{a,b\})$. Consider the first time another move is played:
	
	\begin{itemize}
		\item If the opponent plays an edge in $G' \setminus P'$ containing $a$ or $b$, consider that he played the same move in  $G \setminus P$
		\item If the opponent plays the $(au')$ or $(bv')$, consider that he played , $(au)$ or $(bv)$ respectively.
		\item If the opponent plays an edge in $P'$, let $o, o'$ be the outcome of the subtrees rooted in $u'$ and $v'$, and consider that he plays an edge in $P$ such that the subtrees rooted in $u$ and $v$ in $G$ have the same outcome 
	\end{itemize}
	
	\begin{itemize}
		\item If the optimal move is to play an edge in $G \setminus P$ containing $a$ or $b$, play the same move in  $G' \setminus P'$
		\item If the optimal move is to play $(au)$ or $(bv)$, play $(au')$ or $(bv')$ respectively.
		\item If the optimal move is to play an edge in $P$, let $o, o'$ be the outcome of the subtrees rooted in $u$ and $v$, and play an edge in $P'$ such that the subtrees rooted in $u'$ and $v'$ in $G'$ have the same outcome. 
	\end{itemize}
	
	After these moves, all the remaining vertices of $P$ and $P'$ are trees attached to vertices of $C$. Thus, using Reduction rule~\ref{red: move tree}, we can remove moves in them such that at most one move is playable in each of them (and one remains in $P$ if and only if one remains in $P'$). Thus, we can now apply Reduction rule~ \ref{red: replacing trees} to replace both of them by the same tree, which results into the same graph. Overall, the resulting position if equivalent to the same graph, say $G''$. Thus, as the strategy was a winning strategy for the first (second resp.) played in $G$, it is also one in $G'$ 
\end{proof}

\begin{proof}[Proof of \Cref{thm: ndarck parameterized fen}]
	Let $\alpha$ be the maximum size of the trees in $\mathcal P$. By exhausting application of Reduction rule~\ref{reduction paths}, the resulting graph has at most $(2k-1)$ vertices of degree at least $3$ in $C$, all of them are connected through paths of lengths at most $\alpha$, and each of these vertices is attached to trees of order at most $10$. Overall, after applications of these rules, there are at most $40 \alpha k^2 +\mathcal{O}(k)$ vertices in the resulting graph which has the same outcome as $G$. Thus, it is a quadratic kernel. Moreover, since all the the reduction rules can be computed in quadratic time, we can compute this kernel in quadratic time and compute the winner in time $\mathcal{O}(|E|^2 f(k))$.
\end{proof}

We now study more structured graph classes to obtain efficient algorithms: we provide polynomial-time algorithms for clique trees and a subclass of threshold graphs.

A clique tree is a tree where every vertex is expanded into a clique (which can be of order~1).

\begin{theorem}
	\label{thm-cliquePath}
	\ndarck can be solved in polynomial-time on clique trees.
\end{theorem}

\begin{proof}
	There are three types of vertices in a clique tree: \emph{articulation points} (vertices that, when removed, disconnect the graph: they are the intersection of cliques), \emph{leaves} (vertices of degree~1) and \emph{clique vertices} (the remaining vertices). It is easy to see that the articulation points and leaves form a backbone tree-like structure that behaves like a tree for \ndarck. Such vertices can only be removed if they do not disconnect the graph. For a leaf, this is when their neighbor has only clique vertices as other neighbors (and they get removed together). An articulation point can be removed with a clique vertex if it does not disconnect the graph, or along with another articulation point or leaf vertex under a similar condition. The clique vertices can be removed at any time without any consequence, and indeed they will have to be removed during the course of the game. Hence, the clique vertices are "covering" a tree-like structure, and so there is, as in trees, no strategy to follow: all possible moves will be played, and thus we can easily determine who wins.
\end{proof}

Threshold graphs are a subclass of split graphs, where the vertices can be ordered as $v_1,\ldots,v_n$ such that $N(v_1) \subseteq \ldots \subseteq N(v_n)$, where $N(u)$ denotes the neighborhood of $u$. Alternatively, threshold graphs are graphs that can be constructed by repeatedly adding isolated or universal vertices. Note that this implies that the vertices can be partitioned into a clique and an independent set, meaning threshold graphs are split graphs. (Threshold graphs are actually exactly the intersection of split graphs and cographs, \textit{i.e.}, graphs that contain no induced $P_4$.)  In particular, threshold graphs can be partitioned into cliques $K_1,\ldots,K_m$ and independent sets $S_1,\ldots,S_m$ such that all the vertices of $S_i$ (resp. $K_i$) are twins (\emph{i.e.}, vertices with the exact same neighborhood), $\cup_{i=1}^m K_i$ is a clique, and the vertices of $S_i$ are adjacent to the vertices of $\cup_{j=1}^i K_i$. The construction order is the following: for $i$ from $n$ to 1, add the vertices of $S_i$ as isolated vertices, then the vertices of $K_i$ as universal vertices.
We denote a given threshold graph by $(K \cup S,E)$, where $E$ is the edge set, $K$ is the union of the cliques, and $S$ is the union of the independent sets. If a vertex can be in both $K$ and $S$, we consider that it is in $S$ (so the independent set is \emph{maximal} and the clique is \emph{minimal} considering the fact that they form a partition of the vertices). An example of a threshold graph is depicted in \Cref{fig-threshold}.

\begin{figure}
	\centering
	\begin{tikzpicture}
		\node[noeud] (k0) at (0,0) {};
		\node[noeud] (k1) at (2,0) {};
		\node[noeud] (k2) at (4,0) {};
		\node[noeud] (s0) at (0,2) {};
		\node[noeud] (s1a) at (1.35,2) {};
		\node[noeud] (s1b) at (2.65,2) {};
		\node[noeud] (s2) at (4,2) {};
		\draw (k0)to(s0);
		\draw (k0)to(s1a);
		\draw (k0)to(s1b);
		\draw (k0)to(s2);
		\draw (k1)to(s1a);
		\draw (k1)to(s1b);
		\draw (k1)to(s2);
		\draw (k2)to(s2);
		\draw[bend right = 18] (k0)to(k2);
		\draw (k0)to(k1)to(k2);
		\draw[rounded corners,line width=0.5mm] (-0.5,-0.5) rectangle (4.5,0.5);
		\draw[rounded corners,line width=0.5mm] (-0.5,1.5) rectangle (4.5,2.5);
		\draw[rounded corners,dashed] (1.1,1.75) rectangle (2.9,2.25);
		\draw (-1,0) node[scale=1.5] {$K$};
		\draw (-1,2) node[scale=1.5] {$S$};
	\end{tikzpicture}
	\caption{An example of a threshold graph with a twin-free clique.}
	\label{fig-threshold}
\end{figure}

Another characterization of threshold graphs is as the intersection of split graphs and cographs. Recall that \ndarck is PSPACE-complete on split graphs (\Cref{thm-csgPSPACEcompleteSPLIT}), so the following result implies that there is a complexity gap between clique-twin-free threshold graphs and split graphs:

\begin{theorem}
	\label{thm-thresholdTwinFreeClique}
	Let $T=(K \cup S,E)$ be a threshold graph with $|K|=n$ and such that $K$ is twin-free. If $n \ge 5$, then \ndarck on $T$ is outcome-equivalent to \subgameset{1,2} on a heap of size $n-1$, and thus its outcome can be computed in polynomial time.
\end{theorem}

\begin{proof}
	First, note that, as $S$ is a stable set, it is not possible to remove two vertices of $S$ in a single move. Therefore, the two possible moves are :
	
	\begin{itemize}
		\item Type~1: Removing one vertex of $S$ and an adjacent vertex of $K$.
		\item Type~2: Removing two vertices of $K$.
	\end{itemize}
	
	We first deal with threshold graphs with small cliques. Note that as the clique is supposed to be minimal, any vertex $v \in K$ has at least one neighbor in $S$, otherwise, $v$ could be moved to $S$, and $S$ would stay a stable, but $K$ would be smaller.
	\begin{itemize}
		\item If $T$ is a threshold graph with a twin-free minimum clique of order $1$, then $T$ is a star. Either $T$ is just $P_2$ or $P_3$ (the path on two or three vertices), and the first player wins, or $T$ has at least three edges, and no move is allowed on $T$ thus the second player wins.
		\item If $T$ is a threshold graph with a twin-free minimum clique of order $2$, then, as we have at least two vertices in $S$, type~2 moves are impossible. Moreover, as $K$ is twin-free, and must have a universal vertex, the only possible move is of type~1. After this move, the graph is a star, and we are left with the previous case.
		\item If $T$ is a threshold graph with a twin-free minimum clique of order $3$, then the first player wins by playing a type~2 move, leaving the second player with a star or order at least $4$.
		\item If $T$ is a threshold graph with a twin-free minimum clique of order $4$, then the second player wins. If the first player plays a move of type~$i$, by playing move of type~$(3-i)$, the graph is left with one universal vertex in the clique and at least $3$ leaves connected to it. Therefore, no moves are available, and the second player has won.
	\end{itemize}
	
	To prove the result for $n \ge 5$, we prove that the game will be finished once the clique is reduced to a single vertex.
	
	Let $G = (K\cup S, E)$ be a threshold graph with a maximal independent set $S$ and the clique $K$ of order $n$. Denote by $v_1, \dots, v_n$ the vertices of $K$ such that $N(v_1) \supseteq N(v_2) \supseteq \dots \supseteq N(v_n)$. Similarly, denote by $u_1, \dots, u_N$, with $N \ge n$ the vertices of $S$ such that $N(u_1) \subseteq N(u_2) \subseteq \dots \subseteq N(u_N)$.
	
	Suppose that some moves have already been played. As $|S| \ge |K|$ before the beginning of the game, and as it is not possible to remove more vertices from $S$ than from $K$, at any moment of the game, there will be at least as many vertices remaining in $S$ than in $K$. 
	
	Let $i$ be the largest integer such that $v_i$ has not been removed. By inclusion of the neighborhoods, and as $K$ is a clique and the game is non-disconnecting, $v_i$ is a universal vertex. Therefore, any move that does not remove $v_i$ cannot disconnect the graph, and is legal. Therefore, while $K$ contains at least three vertices, the two types of moves are possible.
	
	Finally, note that, even if $|S| = |K|$ at the beginning of the game, as the optimal move for the player having a winning strategy in \subgameset{1,2}, is to remove the complementary (i.e. playing a move of type~$(3-i)$ against a move of type $i$) to its opponent after each move, in an optimal strategy, one of the three first moves will be a move in the clique, and therefore we will reach a position with $|S| \ge |K| + 2$, and thus the last vertices of the clique will not be removed.  
\end{proof}

\section{Conclusion}
\label{sec-conclusion}

In this paper, we expanded previous PSPACE-hardness results for connected subtraction games, a very large class of vertex deletion games. In particular, we proved that, provided~1 is not in the subtraction set, then the game is PSPACE-complete even on bipartite graphs of any given girth. We also proved that \csgS{k} is PSPACE-complete on split graphs, a very restricted class. However, the cases where~1 is in a subtraction set of at least two elements are still open, as well as other strong conditions on the graph (such as planarity, interval graphs, etc).

We also showed that a natural sufficient winning condition for the second player for \arck is GI-hard, which is one of the few hardness results on this game. While this is not enough to determine the computational complexity of \arck, such an approach might help to understand more about it. A future research direction would be to study more necessary and sufficient conditions for first and second player winning, and determine their complexity.

On the polynomial-time side, we showed that \ndarck is easy on tree-like graphs, such as clique trees, and gave a polynomial kernel parameterized by feedback edge number. Note that this contrasts with \arck, where even subdivided stars with three paths are open. Future research could consider other tree-like classes with unbounded feedback edge number, such as cacti. Furthermore, our study of clique-twin-free threshold graphs seems to imply that the complexity gap may lie between threshold and split graphs. Another direction would be to look at cographs, since threshold graphs are the intersection of split (for which the game is PSPACE-complete) and cographs.

%\clearpage
%References
\bibliographystyle{plain}
\bibliography{biblio}

\begin{thebibliography}{10}

\bibitem{albert2019lessons}
Michael Albert, Richard Nowakowski, and David Wolfe.
\newblock {\em Lessons in play: an introduction to combinatorial game theory}.
\newblock AK Peters/CRC Press, 2019.

\bibitem{beaudou2018octal}
Laurent Beaudou, Pierre Coupechoux, Antoine Dailly, Sylvain Gravier, Julien
  Moncel, Aline Parreau, and Eric Sopena.
\newblock Octal games on graphs: The game 0.33 on subdivided stars and bistars.
\newblock {\em Theoretical Computer Science}, 746:19--35, 2018.

\bibitem{bensmail2022largest}
Julien Bensmail, Foivos Fioravantes, Fionn Mc~Inerney, and Nicolas Nisse.
\newblock The largest connected subgraph game.
\newblock {\em Algorithmica}, 84(9):2533--2555, 2022.

\bibitem{berlekamp2004winning}
Elwyn~R Berlekamp, John~H Conway, and Richard~K Guy.
\newblock {\em Winning ways for your mathematical plays}.
\newblock AK Peters/CRC Press, 2001--2004.

\bibitem{burke2026complexity}
Kyle Burke, Antoine Dailly, and Nacim Oijid.
\newblock Complexity and algorithms for {Arc-Kayles} and {Non-Disconnecting
  Arc-Kayles}.
\newblock In {\em International Conference and Workshops on Algorithms and
  Computation}, 2026.

\bibitem{DBLP:journals/corr/abs-1101-1507}
Kyle Burke and O~George.
\newblock A {PSPACE}-complete graph nim.
\newblock {\em Games of No Chance}, 5:259--269, 2019.

\bibitem{conway2000numbers}
John~H Conway.
\newblock {\em On numbers and games}.
\newblock CRC Press, 2000.

\bibitem{dailly2018criticalite}
Antoine Dailly.
\newblock {\em Criticalit{\'e}, identification et jeux de suppression de
  sommets dans les graphes: Des {\'e}toiles plein les jeux}.
\newblock PhD thesis, Universit{\'e} de Lyon, 2018.

\bibitem{dailly2019generalization}
Antoine Dailly, Valentin Gledel, and Marc Heinrich.
\newblock A generalization of {A}rc-{K}ayles.
\newblock {\em International Journal of Game Theory}, 48:491--511, 2019.

\bibitem{dailly2019connected}
Antoine Dailly, Julien Moncel, and Aline Parreau.
\newblock Connected subtraction games on subdivided stars.
\newblock {\em Integers: Electronic Journal of Combinatorial Number Theory},
  19:G3, 2019.

\bibitem{duchene2024bipartite}
Eric Duch{\^e}ne, Nacim Oijid, and Aline Parreau.
\newblock Bipartite instances of {INFLUENCE}.
\newblock {\em Theoretical Computer Science}, 982:114274, 2024.

\bibitem{DBLP:journals/tcs/Fukuyama03}
Masahiko Fukuyama.
\newblock A {N}im game played on graphs.
\newblock {\em Theoretical Computer Science}, 304(1-3):387--399, 2003.

\bibitem{DBLP:journals/tcs/Fukuyama03a}
Masahiko Fukuyama.
\newblock A {N}im game played on graphs {II}.
\newblock {\em Theoretical Computer Science}, 304(1-3):401--419, 2003.

\bibitem{gardner1974mathematical}
Martin Gardner.
\newblock Mathematical games: Cram, crosscram and quadraphage: new games having
  elusive winning strategies.
\newblock {\em Scientific American}, 230(2):106--108, 1974.

\bibitem{grundy1939mathematics}
Patrick~M Grundy.
\newblock Mathematics and games.
\newblock {\em Eureka}, 2:6--8, 1939.

\bibitem{guy1956g}
Richard~K Guy and Cedric~AB Smith.
\newblock The {G}-values of various games.
\newblock In {\em Mathematical Proceedings of the Cambridge Philosophical
  Society}, volume 52(3), pages 514--526. Cambridge University Press, 1956.

\bibitem{hanaka2024faster}
Tesshu Hanaka, Hironori Kiya, Michael Lampis, Hirotaka Ono, and Kanae
  Yoshiwatari.
\newblock Faster winner determination algorithms for ({C}olored) {A}rc
  {K}ayles.
\newblock {\em Journal of Computer and System Sciences}, page 103716, 2025.

\bibitem{hanaka2023winner}
Tesshu Hanaka, Hironori Kiya, Hirotaka Ono, and Kanae Yoshiwatari.
\newblock Winner determination algorithms for graph games with matching
  structures.
\newblock {\em Algorithmica}, pages 1--17, 2023.

\bibitem{huggan2016polynomial}
Melissa~A Huggan and Brett Stevens.
\newblock Polynomial time graph families for {A}rc {K}ayles.
\newblock {\em Integers}, 16:A86, 2016.

\bibitem{lampis2014computational}
Michael Lampis and Valia Mitsou.
\newblock The computational complexity of the game of set and its theoretical
  applications.
\newblock In {\em LATIN 2014: Theoretical Informatics: 11th Latin American
  Symposium, Montevideo, Uruguay, March 31--April 4, 2014. Proceedings 11},
  pages 24--34. Springer, 2014.

\bibitem{lemoinecomputation}
Julien Lemoine and Simon Viennot.
\newblock Computation records of normal and mis{\`e}re cram.

\bibitem{schaefer1978complexity}
Thomas~J Schaefer.
\newblock On the complexity of some two-person perfect-information games.
\newblock {\em Journal of Computer and System Sciences}, 16(2):185--225, 1978.

\bibitem{siegel2013combinatorial}
Aaron~N Siegel.
\newblock {\em Combinatorial game theory}, volume 146.
\newblock American Mathematical Soc., 2013.

\bibitem{sprague1935mathematische}
Richard Sprague.
\newblock {\"U}ber mathematische kampfspiele.
\newblock {\em Tohoku Mathematical Journal, First Series}, 41:438--444, 1935.

\bibitem{uiterwijk2018construction}
Jos~WHM Uiterwijk.
\newblock Construction and investigation of cram endgame databases.
\newblock {\em ICGA Journal}, 40(4):425--437, 2018.

\bibitem{uiterwijk2019solving}
Jos~WHM Uiterwijk.
\newblock Solving cram using combinatorial game theory.
\newblock In {\em Advances in Computer Games}, pages 91--105. Springer, 2019.

\end{thebibliography}

\end{document}